\documentclass[dvips,ss,preprint]{imsart}

\newtheorem{theorem}{Theorem}[section]
\newtheorem{condition}{Condition}[section]

\newtheorem{proposition}{Proposition}[section]
\begin{document}

\begin{frontmatter}
\title{Regression for partially observed variables and
nonparametric quantiles of conditional probabilities}
\runtitle{Regression for partially observed variables}

\begin{aug}
\author{\fnms{Odile} \snm{Pons}\ead[label=e1]{Odile.Pons@jouy.inra.fr}}
\address{INRA, Math\'ematiques, \\
78352 Jouy en Josas cedex, France\\
\printead{e1}} \runauthor{O. Pons}
\end{aug}

\begin{abstract}
Efficient estimation under bias sampling, censoring or truncation is
a difficult question which has been partially answered and the usual
estimators are  not always consistent. Several biased  designs are
considered for models with variables $(X,Y)$ where $Y$ is an
indicator and $X$ an explanatory variable, or for continuous
variables $(X,Y)$. The identifiability of the models are discussed.
New nonparametric estimators of the regression functions and
conditional quantiles are proposed.
\end{abstract}

\begin{keyword}[class=AMS]
\kwd[Primary ]{62E10} \kwd{62G20} \kwd{62F10}
\end{keyword}

\begin{keyword}
\kwd{Censoring} \kwd{Nonparametric regression} \kwd{Quantile}
\kwd{Truncation}
\end{keyword}
\tableofcontents
\end{frontmatter}

\section{Introduction}\label{introduction}
 Let $(X_i,Y_i)_{i\leq n}$ be a sample of the variable set $(X,Y)$ where
$Y$ is an indicator variable and $X$ is an explanatory variable.
Conditionally on $X$, $Y$ follows a Bernoulli distribution with
parameter $p(x)= \Pr(Y=1| X =x)$. Usual examples are response
variables $Y$ to a dose $X$ or to an expository time $X$, economic
indicators. The variable $X$ may be observed at fixed values $x_i$,
$i\in \{1,\ldots, m\}$ on a regular grid  $\{1/m, \ldots,1\}$ or at
irregular fixed or random times $t_j, j\leq n$, for a continuous
process $(X_t)_{t\leq T}$.

Exponential linear models with known link functions are often used,
especially the logistic regression model defined by
$p(x)=e^{\psi(x)} \{1+e^{\psi(x)}\}^{-1}$ with a parametric function
$\psi$. The inverse function of $p$ is easily estimated using
maximum likelihood estimators of the parameters and many authors
have studied confidence sets for the parameters and the quantiles of
the model.

In a nonparametric setting and for discrete sampling design with
several independent observations for each
value $x_j$ of $X$, the likelihood is written
$$L_n=\prod_{i=1}^n p(X_i)^{Y_i} \{1-p(X_i)\}^{1-Y_i}
= \prod_{j=1}^m \prod_{i=1}^n [\{p(x_j)\}^{Y_i}
\{1-p(x_j)\}^{1-Y_i}]^{ 1_{\{X_i=x_j\}}}.$$
The maximum likelihood estimator of $p(x_j)$
is the proportion of individuals with $Y_i=1$ as $X_i=x_j$,
$$\widehat{p}_{1n}(x_j)= \bar{Y}_{n; x_j}= \frac{1}{\sum_{i=1}^{n}
  1_{\{X_i=x_j \}}}
\sum_{i=1}^{n} Y_i 1_{\{X_i= x_j\}}, j=1,\ldots ,m.$$

Regular versions of this estimator are obtained by kernel
smoothing or by projections on a regular basis of functions,
especially if the variable $X$ is continuous.
Let  $K$ denote a symmetric positive kernel with integral 1,
$h=h_n$ a bandwidth and $K_h(x)= h^{-1} K(h^{-1}x)$, with
$h_n\rightarrow 0$ as $n\rightarrow\infty$.
A local maximum likelihood estimator of $p$ is defined as
$$\widehat{p}_{2n}(x)=\frac{1}{\sum_{i=1}^{n} K_h(x-X_i)} \sum_{i=1}^{n}
Y_i K_h(x-X_i)$$ or by higher order polynomial approximations
\cite{r2}.

Under regularity conditions of $p$ and $K$ and ergodicity of the
process $(X_t, Y_t)_{t\geq 0}$, the estimator $\widehat{p}_{2n}$ is
$P$-uniformly consistent and  asymptotically Gaussian. When $p$ is
monotone, the estimators are asymptotically monotone in probability.
For large $n$, the inverse function $q$ is then estimated by
$\widehat{q}_n(u)= \sup\{x: \widehat{p}_n(x)\leq u\}$ if $p$ is
decreasing or by $\widehat{q}_n(u)= \inf\{x: \widehat{p}_n(x)\geq
u\}$ if $p$ is increasing. The estimator $\widehat{q}_n$ is also
$P$-uniformly consistent and asymptotically Gaussian \cite{r7}. For
small samples, a monotone version of $\widehat{p}_n$ using the
greatest convex minorant or the smallest concave majorant algorithm
may be used before defining a direct inverse. Other nonparametric
inverse functions have been defined \cite{r1}.

Under bias sampling, censoring or truncation, the distribution
function of $Y$ conditionally on $X$ is not always identifiable. The
paper studies several cases and defines new estimators of
conditional and marginal distributions, for a continuous bivariate
set $(X,Y)$ and for a conditional Bernoulli variable $Y$.

\section{Bias depending on the value of $Y$}\label{binomial}
In case-control studies, individuals are not uniformly sampled in
the population: for rare events,  they are sampled  so that the
cases of interest (individuals with $Y_i=1$) are sufficiently
represented in the sample but the proportion of cases in the sample
differs from its proportion in the general population \cite{r6}. Let
$S_i$ be the sampling indicator  of individual $i$ in the global
population and
 $$\Pr(S_i=1| Y_i=1)= \lambda_1, \; \Pr(S_i=1| Y_i=0)= \lambda_0.$$
The distribution function  of $(S_i, Y_i)$ conditionally on $X_i=x$ is given by
\begin{eqnarray*}
\Pr(S_i=1, Y_i=1| x)&=& \Pr(S_i=1| Y_i=1) \Pr(Y_i=1|x)
 =\lambda_1 p(x),  \\
\Pr(S_i=1, Y_i=0| x)&=& \Pr(S_i=1| Y_i=0) \Pr(Y_i=0|x)
 =\lambda_0 \{1-p(x)\}, \\
\Pr(S_i=1|x)&=& \Pr(S_i=1, Y_i=1| x)+\Pr(S_i=1, Y_i=0| x)\\
&=&\lambda_1 p(x)+\lambda_0 \{1-p(x)\}.
\end{eqnarray*}
Let
$$ \theta=\frac{\lambda_0}{ \lambda_1},\;
\alpha(x)=\theta \frac{1-p(x)}{p(x)}.$$
For individual $i$, $(X_i,Y_i)$ is observed
conditionally on $S_i=1$ and the conditional distribution function  of $Y_i$ is
defined by
\begin{eqnarray*}
 \pi(x)&=&\Pr(Y_i=1| S_i=1, X=x)= \frac{\lambda_1 p(x)}{\lambda_1 p(x)+
\lambda_0 \{1-p(x)\}}\\
&=& \frac{p(x)}{p(x)+ \theta
\{1-p(x)\}}= \frac{1}{1+\alpha(x)}.
\end{eqnarray*}
The probability $p(x)$ is deduced from $\theta$ and $\pi(x)$ by the
relation
$$p(x)=\frac{\theta \pi(x)}{1+(\theta-1)\pi(x)}$$ and the
bias sampling is
$$ \pi(x)-p(x)= \frac{(1-\theta)\pi(x)(1-
  \pi(x))}{1+(\theta-1)\pi(x)}.$$

The model defined by $(\lambda_0, \lambda_1, p(x))$ is
over-parameterized and only the function $\alpha$ is  identifiable.
The proportion $\theta$ must therefore be known or estimated from a
preliminary study before an estimation of the probability function
$p$. In  the logistic regression model, $\psi(x)=\log[p(x)  \{1-
p(x)\}^{-1}]$ is replaced by
 $\log \alpha(x)=\log[\pi(x)  \{1-\pi(x)\}^{-1}] =\psi(x)-\log \theta$.
Obviously, the bias sampling modifies the parameters of the model but
not this model and the only  stable parametric model is the logistic
regression.

Let $\gamma$ be the inverse of the proportion of cases in the
population,
\begin{equation}
 \gamma=\Pr(Y=0)/\Pr(Y=1)=E(1-Y)/EY= \frac{1-\int p(x)\, dF_X(x)}{\int
   p(x)\, dF_X(x)}.\label{gamma}
\end{equation}
Under the bias sampling,
\begin{eqnarray*}
 \Pr(Y_i=1| S_i=1)&=&\frac{\lambda_1 \int p\, dF_X}{\lambda_0 (1-\int
   p\, dF_X) +\lambda_1 \int p\, dF_X}=\frac{1}{1+\theta\gamma},\\
 \Pr(Y_i=0| S_i=1)&=&\frac{\lambda_0 (1-\int p\, dF_X)}{\lambda_0 (1-\int
   p\, dF_X )+\lambda_1 \int p\, dF_X}=\frac{\theta\gamma}{1+\theta\gamma},
\end{eqnarray*}
$\gamma$ is modified by the scale parameter $\eta$: it becomes
$\Pr(Y=0| S=1)/\Pr(Y=1| S=1)= \theta\gamma$.

The product $\theta\gamma$ may be directly estimated by maximization
of the likelihood and
$$\theta \widehat{\gamma}_n= 1-\frac{\sum_i Y_i
1_{\{S_i=1\}}}{\sum_i 1_{\{S_i=1\}}}.$$ In a discrete sampling
design with several independent observations for fixed values $x_j$
of the variable $X$, the likelihood is
$$\prod_{i=1}^n \pi(X_i)^{Y_i} \{1-\pi(X_i)\}^{1-Y_i}=
\prod_{j=1}^m\prod_{i=1}^n [\pi(X_i)^{Y_i} \{1-\pi(X_i)\}^{1-Y_i}]^{
1_{\{X_i=x_j \}}}$$
and $\alpha_j=\alpha(x_j)$ is estimated by
$$\widehat{\alpha}_{jn}= \frac{\sum_i (1-Y_i) 1_{\{S_i=1\}}
  1_{\{X_i=x_j \}}}{\sum_i Y_i 1_{\{X_i= x_j\}} 1_{\{S_i=1\}}}.$$
For random observations of the variable $X$, or for fixed observations
without replications, $\alpha(x)$ is estimated by
$$\widehat{\alpha}_{n}(x)=\frac{\sum_i (1-Y_i) 1_{\{S_i=1\}} K_h(x-X_i)}{\sum_i
Y_i 1_{\{S_i=1\}} K_h(x-X_i)}.$$
If $\theta$ is known, nonparametric estimators of $p$ are deduced as
\begin{eqnarray*}
\widehat{p}_{n}(x_j)&=& \frac{\theta \sum_i Y_i 1_{\{S_i=1\}} 1_{\{X_i=x_j \}}}
{\sum_i (1-Y_i+\theta Y_i) 1_{\{S_i=1\}} 1_{\{X_i=x_j \}}}, \mbox{ in
  the discrete case},\\
\widehat{p}_{n}(x)&=& \frac{\theta\sum_i Y_i 1_{\{S_i=1\}} K_h(x-X_i)}{\sum_i
(1-Y_i+\theta Y_i) 1_{\{S_i=1\}} K_h(x-X_i)}, \mbox{ in
  the  continuous case}.
\end{eqnarray*}

\section{Bias due to truncation on $X$}\label{Xtruncation}
Consider that $Y$ is observed under a fixed
truncation of $X$: we assume that $(X, Y)$ is observed only
if $X\in [a,b]$, a sub-interval of the support $I_X$ of
the variable $X$, and $S=1_{[a,b]}(X)$. Then
$$\Pr(Y_i=1)= \int_{I_X} p(x)\, dF_X(x),\quad
\Pr(Y_i=1, S_i=1)= \int_{a}^b p(x)\, dF_X(x)$$
and the conditional probabilities of sampling, given the status value, are
\begin{eqnarray*}
\lambda_1&=&\Pr(S_i=1|Y_i=1)= \frac{\int_{a}^b p(x)\, dF_X(x)}{\int_{I_X}
p(x)\, dF_X(x)},  \\
\lambda_0&=&\Pr(S_i=1|Y_i=0)= \frac{\int_{a}^b \{1- p(x)\} \,
  dF_X(x)}{1-\int_{I_X} p(x) \, dF_X(x)}.
\end{eqnarray*}
If the ratio $\theta=\lambda_0/\lambda_1$ is known or  otherwise
estimated, the previous estimators may be
used for the estimation of $p(x)$ from the truncated sample with
$S_i\equiv 1$.

For a random truncation interval $[A, B]$, the sampling indicator is
$S=1_{[A, B]}(X)$ and the integrals of $p$ are replaced by their
expectation with respect to the distribution function  of $A$ and $B$ and the
estimation is similar.\\

\section{Truncation of a response variable in a nonparametric regression
model} \label{Ytruncation}
Consider then $(X,Y)$ a two-dimensional variable in a
left-truncated transformation model:
Let  $Y$ denote a response to a continuous expository variable $X$, up to a
variable of individual variations $\varepsilon$ independent of $X$,
$$ Y=m(X)+\varepsilon,\quad E\varepsilon =0, \quad E
\varepsilon^{2}<\infty,$$
$(X, \varepsilon)$ with distribution function  $(F_X, F_{\varepsilon})$.
The distribution function  of $Y$ conditionally on  $X$ is defined~by
\begin{eqnarray}
F_{Y|X}(y; x) &=& P(Y\leq y| X=x)=F_{\varepsilon }(y-m(x)),\label{model} \\
m(x) &=& E(Y|X=x),\nonumber
\end{eqnarray}
and the function $m$ is continuous.
The joint and marginal distribution functions of $X$ and $Y$ are denoted
$F_{X, Y}$, with support $I_{Y,X}$, $F_X$, with bounded support $I_X$,
and $F_Y$, such that $F_{Y}(y)= \int F_\varepsilon(y-m(s))\,  dF_X(s)$ and
$F_{X, Y}(x,y)=  \int 1_{\{ s \leq x \}}
F_\varepsilon(y-m(s))\,  dF_X(s).$

The observation of $Y$ is supposed left-truncated by a variable $T$
independent of $(X, Y)$, with distribution function  $F_T$ : $Y$ and $T$
are observed conditionally on $Y\geq T$ and none of the variables
is observed if  $Y< T$. Denote $\bar{F}=1- F$ for any distribution
function  $F$ and, under left-truncation,
\begin{eqnarray}
\alpha(x)  &=& P(T \leq  Y|X = x) =   \int_{-\infty}^{\infty}
\bar{F}_\varepsilon(y-m(x)) \, dF_T(y),
\nonumber \\
A(y;x)  &=& P(Y\leq y| X=x, T \leq  Y) \nonumber \\
&=& \alpha^{-1}(x) \int_{-\infty}^{y}  F_T(v) \,
dF_\varepsilon(v-m(x)) \label{Ayx} \\
B(y;x) &=& P(T \leq y \leq  Y| X=x, T \leq  Y) \nonumber \\
&=& \alpha^{-1}(x) F_T(y) \bar{F}_\varepsilon(y -m(x)), \label{Byx} \\
m^*(x) &=& E(Y| X=x, T \leq  Y) =
 \alpha^{-1}(x) \int y   F_T(y) \,dF_{Y|X}(y;x). \nonumber
\end{eqnarray}
Obviously, the mean of $Y$ is biased under the truncation and a
direct estimation of the conditional distribution function $F_{Y|X}$
is of interest for the estimation of $m(x)=E(Y|X=~x)$ instead of the
apparent mean $m^*(x)$. The function $\bar{F}_\varepsilon$ is also
written $\exp \{- \Lambda_\varepsilon\}$ with $
\Lambda_\varepsilon(y) = \int_{-\infty}^{y} \bar{F}_\varepsilon^{-1}
dF_\varepsilon$ and the expressions (\ref{Ayx})-(\ref{Byx}) of $A$
and $B$ imply that
$$\Lambda_\varepsilon(y - m(x)) =  \int_{-\infty}^{y} B^{-1}(s;x)
 \,A(ds; x)$$
and  $F_{Y|X}(y; x) = \exp \{- \Lambda_\varepsilon(y-m(x))\}$.

An estimator of $F_{Y|X}(y; x)$ is obtained as the product-limit estimator
$\widehat{F}_{\varepsilon,n}(y- m(x))$ of
$F_\varepsilon(y-m(x))$ based on estimators of $A$ and $B$: For a
sample $(X_i, Y_i)_{1\leq i\leq n}$, let $x$ in
$I_{X,n,h}=[\min_i X_i+h, \max_i X_i-h]$ and
\begin{eqnarray}
 \widehat{A}_n(y;x) &=& \frac{\sum_{i=1}^n
K_{h}(x-X_i) I_{\{T_i\leq Y_i\leq y\}}}
{\sum_{i=1}^n K_{h}(x-X_i) I_{\{T_i\leq Y_i\}}}, \nonumber\\
\widehat{B}_n(y; x) &=&  \frac{\sum_{i=1}^n
K_{h}(x-X_i) I_{\{T_i\leq y \leq Y_i\}}}
{\sum_{i=1}^n K_{h}(x-X_i) I_{\{T_i\leq Y_i\}}}; \nonumber\\
\widehat{F}_{Y|X,n}(y;x) &=&1- \prod_{1\leq Y_i\leq y}
\left\{1 - \frac{d\widehat{A}_{n}}{\widehat{B}_n}(Y_i; x)
\right\} \nonumber\\
&=& 1- \prod_{1\leq i\leq n} \left\{1- \frac{ K_{h}(x-X_i)
  I_{\{T_i\leq Y_i\leq y\}}}{\sum_{j=1}^n
K_{h}(x-X_j) I_{\{T_j\leq Y_i \leq Y_j\}}} \right\},\label{estim_Fyx}
\end{eqnarray}
with $0/0=0$.
That is a nonparametric maximum likelihood estimator of $F_{Y|X}$, as is the
Kaplan-Meier estimator for the distribution function  of a
right-censored variable.
Then an estimator of $m(x)$ may be defined as an estimator of $\int y
\, F_{Y|X}(dy;x)$,
\begin{eqnarray}
\widehat{m}_n(x) &=&  \sum_{i=1}^n  Y_i  I_{\{T_i\leq Y_i\}}
\{\widehat{F}_{Y|X,n}(Y_i; x) -\widehat{F}_{Y|X,n}(Y^-_i; x)\}\nonumber\\
& =&  \frac{\sum_{i=1}^n  Y_i I_{\{T_i\leq Y_i\}}
 K_{h}(x-X_i) \, \widehat{F}_{Y|X,n}(Y^-_i; x)}{\sum_{j=1}^n
  K_{h}(x-X_j) I_{\{T_j\leq Y_i\leq Y_j\}}}.\label{estim_m}
\end{eqnarray}
By the same arguments, from means in (\ref{Ayx})-(\ref{Byx}),
$\bar{F}_Y(y) =
  E\bar{F}_\varepsilon(y-m(X))$ is estimated by
$$\widehat{\bar{F}}_{Y,n}(y) = \prod_{1\leq i\leq n} \left\{1- \frac{
  I_{\{T_i\leq Y_i\leq y\}}}{\sum_{j=1}^n
I_{\{T_j\leq Y_i\leq Y_j\}}} \right\},$$ the distribution function
$F_T$ is simply estimated by the product-limit estimator for
right-truncated variables \cite{r10}
$$\widehat{F}_{T,n}(t)=\prod_{1\leq i\leq n} \left\{1 -\frac{ I_{\{t \leq
    T_i\leq Y_i\}}}{\sum_{j=1}^n I_{\{T_j\leq T_i\leq Y_j\}}}
    \right\}$$
and an estimator of $F_\varepsilon$ is deduced from those of
$F_{Y|X}$, $F_X$ and $m$ as
$$\widehat{F}_{\varepsilon,n}(s)= n^{-1}\sum_{1\leq i\leq n}
\widehat{F}_{Y|X, n}(s+\widehat{m}_n(X_i); X_i).$$
The means of $T$ and $C$ are estimated by
$$\widehat{\mu}_{T,n} =  n^{-1}\sum_{i=1}^n  \frac{T_i  I_{\{T_i\leq Y_i\}}\,
  \widehat{F}_{T,n}(T^-_i)}{\sum_{j=1}^n I_{\{T_j\leq T_i\leq
  Y_j\}}},\,
\widehat{\mu}_{Y,n} =  n^{-1}\sum_{i=1}^n  \frac{Y_i  I_{\{T_i\leq Y_i\}}\,
  \widehat{F}_{Y,n}(Y^-_i)}{\sum_{j=1}^n I_{\{T_j\leq Y_i\leq Y_j\}}}.$$
The estimators $\widehat{F}_{Y,n}$ and $\widehat{F}_{T,n}$ are known
to be $P$-uniformly consistent and asymptotically Gaussian. For the
further convergence restricted to the interval $I_{n,h}=\{(y,x)\in
I_{Y,X}: x\in I_{X,n,h}\}$, assume
\begin{condition}
C1. $h=h_n\rightarrow 0$  and $nh^3\rightarrow \infty$ as
$n\rightarrow\infty$,
 $\int K =1$, $\kappa_1= \int x^2 K(x)\, dx$ and  $\kappa_2= \int K^2< \infty$. \\
C2. the conditional probability $\alpha$ is strictly positive in the
interior of $I_X$,\\
 C3. The distribution function  $F_{Y,X}$ is
twice continuously differentiable with respect to $x$ and
differentiable with respect to $y$.\\
C4. $E\varepsilon^{2+\delta}< \infty$ for a $\delta$ in $(1/2, 1]$.
\end{condition}
Let us denote $\dot{F}_{Y,X,2}(y,x) = \partial F_{Y|X}(y,
x)/\partial x$,   $\ddot{F}_{Y,X,2}(y,x) =
\partial^2 F_{Y|X}(y, x)/\partial x^2$,
and $\dot{F}_{Y|X,1}(y,x)= \partial F_{Y|X}(y, x)/\partial y$.
\begin{proposition}\label{cvAnBn}
 $\sup_{I_{n,h}}|\widehat{A}_n-A|\stackrel{P}{\rightarrow} 0$ and
 $\sup_{I_{n,h}} |\widehat{B}_n-B|\stackrel{P}{\rightarrow} 0$,
\begin{eqnarray*}
 b_{nh}^A(y; x) &\equiv& (E\widehat{A}_n-A)(y; x)=
 \frac{h^2}{2\alpha(x)}  \kappa_1
 \left\{\int_{-\infty}^y F_T(v) \ddot{F}_{Y,X,2}(dv,dx) \right.\\
&& \left. - A(y;x) \int_{-\infty}^\infty F_T(v)
 \ddot{F}_{Y,X,2}(dv,dx)\right\}  +o(h^2),\\
 b_{nh}^B(y; x) &\equiv& (E\widehat{B}_n-B)(y; x)=
\frac{h^2}{2\alpha(x)}  \kappa_1
 \{F_T(y) \int \ddot{F}_{Y,X,2}(dv,dx)\, dx\\
&&- B(y;x) \int  F_T(v) \ddot{F}_{Y,X,2}(dv,dx)\}  +o(h^2),\\
v_{nh}^A(y; x) &\equiv&{\rm var} \widehat{A}_n(y; x)=
(nh)^{-1}  \kappa_2 A(1-A)(y;x) \alpha^{-1}(x) +
 o((nh)^{-1}),\\
v_{nh}^B(y; x) &\equiv&{\rm var} \widehat{B}_n(y; x)=
(nh)^{-1}  \kappa_2 B(1-B)(y;x) \alpha^{-1}(x) +
 o((nh)^{-1}).
\end{eqnarray*}
If $nh^5\rightarrow 0$,
$(nh)^{1/2}(\widehat{A}_n-A)$ and $(nh)^{1/2}(\widehat{B}_n-B)$
 converge in distribution to Gaussian
processes with mean zero, variances  $\kappa_2 A(1-A)(y;x)
\alpha^{-1}(x)$ and $\kappa_2 B(1-B)(y;x) \alpha^{-1}(x)$
respectively, and the covariances of the limiting processes are zero.
\end{proposition}
The proof relies on an expansion of the form
$$(nh)^{1/2}(\widehat{A}_n-A)(y;x)=(nh)^{1/2} c^{-1}(x)
\{(\widehat{a}_n-a)(y;x) - A (\widehat{c}_n- c)(x)\}+o_{L^2}(1)$$
with  $\widehat{A}_n= \widehat{c}_n^{-1} \widehat{a}_n$ and
$\widehat{B}_n= \widehat{c}_n^{-1} \widehat{b}_n$, where
$$\widehat{c}_n(x) = n^{-1}
\sum_{i=1}^n K_{h}(x-X_i) I_{\{T_i\leq Y_i\}},\,
\widehat{a}_n(y;x)=n^{-1} \sum_{i=1}^n
K_{h}(x-X_i) I_{\{T_i\leq Y_i\leq y\}},$$
$$\widehat{b}_n(y; x)=n^{-1} \sum_{i=1}^n
K_{h}(x-X_i) I_{\{T_i\leq y \leq Y_i\}}.$$ A similar approximation
holds for $\widehat{B}_n$. The biases and variances are deduced from
those of each term and the weak convergences are established as in
\cite{r7}.

From proposition \ref{cvAnBn} and applying the results of the
nonparametric regression,
\begin{proposition}\label{cvFnmn}
The estimators
$\widehat{F}_{Y|X, n}$, $\widehat{m}_n$, $\widehat{F}_{\varepsilon,n}$
converge $P$-uniformly to $F_{Y|X}$, $m$, $F_\varepsilon$,
$\widehat{\mu}_{Y, n}$ and $\widehat{\mu}_{T, n}$ converge
$P$-uniformly to $EY$ and $ET$ respectively.
\end{proposition}
The weak convergence of the estimated distribution function  of
truncated survival data was proved in several papers \cite{r4, r5}.
As in \cite{r3} and by proposition \ref{cvAnBn}, their proof extends
to their weak convergence on $(\min_i \{Y_i: T_i<Y_i\}, \max_i
\{Y_i: T_i<Y_i\})$ under the conditions  $\int F_T \, dF_{Y|X} <
\infty$ and $\int \bar{F}_{Y|X}^{-1}\,dF_T < \infty$ on $I_{X,n,h}$,
which are simply satisfied if for every $x$ in $I_{X,n,h}$, $\inf
\{t: F_T(t) >0\} <\inf \{t: F_{Y|X}(t;x) >0\}$ and $\sup \{t:
F_{Y|X}(t;x) >0\} < \sup \{t: F_T(t) >0\}$.

\begin{theorem}\label{cvFyx}
$(nh)^{1/2}(\widehat{F}_{Y|X,n}-  F_{Y|X})1_{I_{n,h}}$ converges
weakly to a centered Gaussian process $W$ on $I_{Y,X}$. The
variables
 $(nh)^{1/2}(\widehat{m}_n- m)(x)$, for every $x$ in $I_{X,n,h}$, and
$(nh)^{1/2}(\widehat{\mu}_{Y, n}- E Y)$
converge weakly to $E W(Y;x)$ and $E \int W(Y;x)\, dF_X(x)$.
\end{theorem}

If $m$ is supposed monotone  with  inverse function $r$, $X$ is
written $X=r(Y-\varepsilon)$ and the quantiles of $X$ are defined by
the inverse functions $q_1$ and $q_2$ of $F_{Y|X}$ at fixed $y$ and $x$,
respectively, are defined by the equivalence between
$$ F_{Y|X}(y; x) = u \quad \mbox{ and } \quad
\left\{\begin{array}{ll}
x&=   r(y-Q_\varepsilon(u))\equiv q_1(u; y)\\
y&= m(x) + Q_\varepsilon(u)\equiv q_2(u;x),
\end{array} \right.$$
where $Q_\varepsilon(u)$ is the inverse of $F_\varepsilon$ at $u$.
Finally, if $m$ is increasing, $F_{Y|X}(y;x)$ is decreasing in $x$ and
increasing in $y$, and it is the same for its estimator
$\widehat{F}_{Y|X,n}$, up to a random set of small probability.
The thresholds $q_1$ and $q_2$ are estimated by
\begin{eqnarray*}
\widehat{q}_{1,n,h}(u;y) &=& \sup \{ x: \widehat{F}_{Y|X,n}(y;x)\leq u \},\\
\widehat{q}_{2,n,h}(u;x) &=&  \inf \{ y: \widehat{F}_{Y|X,n}(u;x)\geq u \}.
\end{eqnarray*}
As a  consequence of theorem \ref{cvFyx} and generalizing known results on
quantiles
\begin{theorem}
For $k=1, 2$, $\widehat{q}_{k,n,h}$  converges $P$-uniformly to
 $q_k$ on $\widehat{F}_{Y,X,n}(I_{n,h})$.
For every $y$ and (respect.) $x$, $(nh)^{1/2}(\widehat{q}_{1,n,h}-
 q_1)(\cdot; y)$ and
$(nh)^{1/2}(\widehat{q}_{2,n,h}-  q_2)(\cdot; x)$ converge weakly to
the centered Gaussian process $W\circ q_1[\dot{F}_{Y|X,1}(y;
q_1(\cdot; y)]^{-1}$ and,
 respect., $W\circ q_2[\dot{F}_{Y|X,2}(q_2(\cdot; x); x)]^{-1}$.
\end{theorem}

\section{Truncation and censoring of $Y$  in a
  nonparametric model}\label{Ytruncationcensoring}
The variable $Y$ is supposed left-truncated by  $T$ and right-censored
by a variable $C$ independent of $(X, Y,T)$.
The notations $\alpha$ and those of the joint and marginal
distribution function  of $X$,  $Y$ and $T$ are in section \ref{Ytruncation}
and $F_C$ is the distribution function  of $C$.   The observations are
$\delta=1_{\{Y\leq C\}}$, and $(Y\wedge C, T)$, conditionally on
$Y\wedge C\geq T$. Let
\begin{eqnarray*}
A(y;x)  &=& P(Y\leq y\wedge C| X=x, T \leq Y) \\
&=& \alpha^{-1}(x) \int_{-\infty}^{y}  F_T(v)  \bar{F}_C(v)\,
F_{Y|X}(dv;x) \\
B(y;x) &=& P(T \leq y \leq  Y\wedge  C| X=x, T \leq  Y)  \\
&=& \alpha^{-1}(x) F_T(y) \bar{F}_C(y) \bar{F}_{Y|X}(y;x),\\
\bar{F}_{Y|X}(y; x) &=& \exp \{- \int_{-\infty}^y B^{-1}(v;x)\, A(dv;x)\}.
\end{eqnarray*}
The estimators are now written
\begin{eqnarray*}
\widehat{\bar{F}}_{Y|X,n}(y;x) &=&
\prod_{1\leq i\leq n} \left\{1- \frac{ K_{h}(x-X_i)
  I_{\{T_i\leq Y_i\leq y\wedge C_i\}}}{\sum_{j=1}^n
K_{h}(x-X_j) I_{\{T_j\leq Y_i \leq Y_j\wedge C_j\}}} \right\},\\
\widehat{m}_n(x) &=& \frac{\sum_{i=1}^n  Y_i I_{\{T_i\leq Y_i\leq C_i\}}
 K_{h}(x-X_i) \, \widehat{F}_{Y|X,n}(Y^-_i; x)}{\sum_{j=1}^n
  K_{h}(x-X_j) I_{\{T_j\leq Y_i\leq Y_j\wedge C_j\}}},\\
\widehat{\bar{F}}_{Y,n}(y) &=& \prod_{1\leq i\leq n} \left\{1- \frac{
  I_{\{T_i\leq Y_i\leq y\wedge C_i\}}}{\sum_{j=1}^n
I_{\{T_j\leq Y_i \leq Y_j\wedge C_j\}}} \right\}.
\end{eqnarray*}

If $Y$ is only right-truncated by
$C$ independent of  $(X, Y)$, with observations $(X, Y)$ and
$C$  conditionally on $Y\leq C$, the expressions $\alpha$, $A$ and $B$
are now written
\begin{eqnarray*}
\alpha(x)  &=& P(Y \leq C|X = x) =\int_{-\infty}^{\infty}
\bar{F}_C(y) \, F_{Y|X}(dy;x),  \\
A(y;x)  &=& P(Y\leq y| X=x, Y \leq C)
=\alpha^{-1}(x) \int_{-\infty}^{y}  \bar{F}_C(v)\,
F_{Y|X}(dv;x), \\
B(y;x) &=& P(Y\leq y \leq  C| X=x, Y \leq C)  =
 \alpha^{-1}(x) \bar{F}_C(y) F_{Y|X}(y;x),\\
A'(y;x) &=& P(Y\leq  C\leq y | X=x, Y \leq C)\\
& =& \alpha^{-1}(x) \int_{-\infty}^{y} F_{Y|X}(v;x)\, dF_C(v).
\end{eqnarray*}
The distribution function  $F_C$ and $F_{Y|X}$ are both identifiable
 and their expression differs from the previous ones,
\begin{eqnarray*}
 \bar{F}_C &=&  \exp \{- \int_{-\infty}^{\cdot}  E B^{-1}(v;X)\,
 E A'(dv;X)\},\\
F_{Y|X}(\cdot;x)& =& \exp \{- \int_\cdot^{\infty}  B^{-1}(v;x)\,
 A(dv;x)\}.
\end{eqnarray*}
The estimators are now
\begin{eqnarray*}
 \widehat{F}_{Y|X,n}(y;x) &=&
\prod_{1\leq i\leq n} \left\{1- \frac{ K_{h}(x-X_i)
  I_{\{Y_i\leq y\wedge C_i\}}}{\sum_{j=1}^n
K_{h}(x-X_j) I_{\{Y_j\leq Y_i\leq C_j\}}} \right\},\\
 \widehat{\bar{F}}_{C,n}(y) &=& \prod_{1\leq i\leq n} \left\{1-
 \frac{I_{\{Y_i\leq C_i\leq y\}}}{\sum_{j=1}^n I_{\{Y_j\leq Y_i\leq
     C_j\}}} \right\},\\
\widehat{m}_n(x) &=& \frac{\sum_{i=1}^n  Y_i I_{\{Y_i\leq C_i\}}
 K_{h}(x-X_i) \, \widehat{F}_{Y|X,n}(Y^-_i; x)}{\sum_{j=1}^n
  K_{h}(x-X_j) I_{\{Y_j\leq Y_i\leq C_j\}}}.
\end{eqnarray*}
If $Y$ is left and right-truncated by variables $T$ and
$C$ independent and independent of  $(X, Y)$, the observations are
$(X,Y)$, $C$ and $T$, conditionally on $T\leq Y\leq C$,
\begin{eqnarray*}
\alpha(x)  &=& P(T \leq  Y \leq C|X = x) =\int_{-\infty}^{\infty}
F_T(y) \bar{F}_C(y) \, F_{Y|X}(dy;x),  \\
A(y;x)  &=& P(Y\leq y| X=x, T \leq Y \leq C) \\
&=& \alpha^{-1}(x) \int_{-\infty}^{y}  F_T(v)  \bar{F}_C(v)\,
F_{Y|X}(dv;x), \\
B(y;x) &=& P(T \leq y \leq  Y| X=x, T \leq  Y \leq C)  \\
&=& \alpha^{-1}(x) F_T(y) \int_y^{\infty} \bar{F}_C(v)\,
\,F_{Y|X}(dv;x),\\
A^\prime(y)  &=& P(y\leq T | T \leq Y \leq C)
= \int_y^{\infty}  dF_T(t)   \int_t^{\infty} \bar{F}_C\, dF_Y, \\
B'(y) &=& P(Y\leq y \leq  C| T \leq Y \leq C)  =
\bar{F}_C(y) \int_{-\infty}^y F_T \, dF_Y, \\
B^{\prime\prime}(y) &=& P(C \leq y| T \leq Y \leq C)
=\int_{-\infty}^y \{\int_{-\infty}^s F_T(v)\, dF_Y(v) \}\, d F_C(s).
\end{eqnarray*}
The distribution functions $F_C$, $F_T$ and  $F_{Y|X}$ are
identifiable, with $F_{Y|X}$ defined by $F_{Y|X}(y;x)=
-\int_{-\infty}^{y} \bar{F}_C^{-1}\, dH(\cdot;x)$ and
\begin{eqnarray*}
H(y;x)&\equiv &\int_y^{\infty} \bar{F}_C(v)\, dF_{Y|X}(dv;x) =\exp \{-
 \int_{-\infty}^{y}  B^{-1}(v;x)\, A(dv;x)\},\\
\bar{F}_C(s)&=& \exp \{- \int \int_{-\infty}^{s}  B^{\prime-1}\,
dB^{\prime\prime}\},\\
F_T(t) &=&\exp [- \{\int_t^\infty  (E B(\cdot;X))^{-1}\, dA^{\prime}\}].
\end{eqnarray*}
Their estimators are
\begin{eqnarray*}
\widehat{\bar{F}}_{C,n}(s)&=&  \prod_{i=1}^n
\left\{1- \frac{ I_{\{T_i\leq Y_i\leq  C_i\leq s\}}}{\sum_{j=1}^n
I_{\{T_j\leq Y_j\leq  C_i\leq C_j\}}} \right\},\\
\widehat{F}_{T,n}(t)&=& \prod_{i=1}^n
\left\{1- \frac{ I_{\{T_i\leq Y_i\leq  C_i\leq t\}}}{\sum_{j=1}^n
  I_{\{T_j\leq T_i\leq   Y_j\leq C_j\}}} \right\},\\
\widehat{F}_{Y|X}(y;x)&=&  \frac{\sum_{i=1}^n
 \widehat{\bar{F}}^{-1}_{C,n}(Y_i)  I_{\{T_i\leq Y_i\leq  C_i\wedge y\}}
 K_{h}(x-X_i) \, \widehat{H}_{Y|X,n}(Y^-_i; x)}{\sum_{j=1}^n
  K_{h}(x-X_j) I_{\{T_j\leq Y_i\leq Y_j\leq C_j\}}};\\
\widehat{H}_{Y|X}(y;x)&=&  \prod_{i=1}^n
\left\{1- \frac{K_{h}(x-X_i) I_{\{T_i\leq Y_i\leq  C_i\wedge
    y\}}}{\sum_{j=1}^n  K_{h}(x-X_j) I_{\{T_j\leq
Y_i \leq Y_j\leq C_j\}}} \right\}.
\end{eqnarray*}
The other nonparametric estimators of the introduction and the
results of section \ref{Ytruncation} generalize to all the
estimators of this section.

Right and left-truncated distribution functions $F_{Y|X}$ and the
truncation distributions are estimated in a closed form by the
solutions a self-consistency equation \cite{r8,r9}. The estimators
still have asymptotically Gaussian limits even with dependent
truncation distributions, when the martingale theory for point
processes does not apply.

\section{Observation by interval}
Consider model (\ref{model}) with an independent censoring variable
$C$ for $Y$. For observations by intervals, only $C$ and the
indicators that $Y$ belongs to the interval $]-\infty, C]$ or $]C,
\infty[$ are observed. The function $F_{Y|X}$ is not directly
identifiable and efficient estimators for $m$ and $F_{Y|X}$ are
maximum likelihood estimators. Let $\delta=I_{\{Y\leq C\}}$ and
assume that $F_\varepsilon$ is $C^2$. Conditionally on $C$ and
$X=x$, the log-likelihood of $(\delta, C)$ is $$l(\delta, C) =
\delta \log F_\varepsilon(C-m(x)) + (1-\delta) \log
\bar{F}_\varepsilon(C-m(x))$$ and its derivatives with respect to
$m(x)$ and $F_\varepsilon$ are
\begin{eqnarray*}
\dot{l}_{m(x)}(\delta, C) &=& - \delta
\frac{f_\varepsilon}{F_\varepsilon}(C-m(x)) +  (1-\delta)
\frac{f_\varepsilon}{\bar{F}_\varepsilon}(C-m(x)),\\
\dot{l}_\varepsilon a(\delta, C) &=& \delta
\frac{\int_{-\infty}^{C-m(x)} a\,
  dF_\varepsilon}{F_\varepsilon(C-m(x))} + (1-\delta)
\frac{\int_{C-m(x)}^\infty a\, dF_\varepsilon}{\bar{F}_\varepsilon(C-m(x))}
\end{eqnarray*}
for every $a$ s.t. $\int a\, dF_\varepsilon=0$ and $\int a^2\,
dF_\varepsilon<\infty$.  With $a_F= -f'_\varepsilon f^{-1}_\varepsilon$,
$\dot{l}_\varepsilon a_F =  \dot{l}_{m(x)}$ then $\dot{l}_{m(x)}$ belongs
to the tangent space for $F_\varepsilon$  and the estimator of $m(x)=
E (Y| X=x)$ must be determined from the estimator of $F_\varepsilon$
through the conditional probability function of the observations
$$B(t;x)= P(Y\leq C\leq t|X=x)=
\int_{-\infty}^t  F_\varepsilon(s-m(x)) \, dF_C(s).$$

Let $\widehat{F}_{C,n}$ the empirical estimator of $F_C$ and
$$\widehat{B}_n(t;x)=  \frac{\sum_{i=1}^n K_{h}(x-X_i) I_{\{Y_i\leq
    C_i\leq t\}}}{\sum_{i=1}^n K_{h}(x-X_i)},$$
an estimator $\widehat{F}_{\varepsilon,n}(t-m(x))$ of
    $F_{\varepsilon,n}(t-m(x))$ is deduced by deconvolution and
$$\widehat{m}_n(x)= \int t\, d\widehat{F}_{\varepsilon,n}(t-m(x)).$$

\end{document}